\documentclass[twoside,11pt,a4paper,leqno]{article}
\usepackage{amsfonts}
\usepackage{amsfonts}
\usepackage{amsfonts}
\usepackage{amsfonts}
\usepackage{amssymb,amsmath,amscd,euscript,verbatim,array}
\usepackage[center,pagestyles]{titlesec}
\usepackage{geometry}
\usepackage{anysize}
\usepackage{fancyhdr}
\usepackage{indentfirst}
\usepackage{graphicx}
\usepackage{color}
\usepackage{ifpdf}
\marginsize{3.5cm}{3.5cm}{2cm}{2cm}

\titleformat{\section}{\centering\normalsize}{\thesection.}{0.5em}{}
\titleformat{\subsection}{\normalsize\bfseries}{\thesubsection.}{0.5em}{}
\titleformat{\subsubsection}{\normalsize\bfseries}{\thesubsubsection.}{0.5em}{}
\newcommand{\N}{\mathbb{N}}
\newcommand{\Z}{\mathbb{Z}}
\newcommand{\R}{\mathbb{R}}
\newcommand{\Q}{\mathbb{Q}}
\newcommand{\C}{\mathbb{C}}

\newtheorem{Theorem}{Theorem}[section]
\newtheorem{Definition}[Theorem]{Definition}
\newtheorem{Lemma}[Theorem]{Lemma}
\newtheorem{Exercise}[Theorem]{Exercise}
\newtheorem{Proposition}[Theorem]{Proposition}
\newtheorem{Remark}[Theorem]{Remark}

\newcommand{\T}{\mathbb{T}}
\newcommand{\bthm}{\begin{Theorem}}
\newcommand{\ethm}{\end{Theorem}}
\newcommand{\bpr}{\begin{Proposition}}
\newcommand{\epr}{\end{Proposition}}
\newcommand{\blm}{\begin{Lemma}}
\newcommand{\elm}{\end{Lemma}}
\newcommand{\bex}{\begin{Exercise}}
\newcommand{\eex}{\end{Exercise}}
\newcommand{\be}{\begin{equation}}
\newcommand{\ee}{\end{equation}}
\newcommand{\beal}{\begin{aligned}}
\newcommand{\enal}{\end{aligned}}
\newcommand{\brm}{\begin{Remark}}
\newcommand{\erm}{\end{Remark}}
\newcounter{item}[section]

\newcommand{\Proof}{\textbf{Proof}\hspace{0.3cm}}
\newcommand{\End}{\ensuremath{\hfill{\Box}}\\}
\renewcommand{\title}[1]{\begin{center}\textbf{\large #1}\end{center}}
\renewcommand{\author}[1]{\begin{center}\small #1\end{center}}
\renewcommand{\date}[1]{\begin{center}#1\end{center}}

\makeatletter \@addtoreset{equation}{section}
\makeatother
\begin{document}
\vspace{10pt}
\title{DESTRUCTION OF INVARIANT CIRCLES FOR GEVREY AREA-PRESERVING TWIST MAPS}
\vspace{6pt}
\author{\sc Lin Wang}
\vspace{10pt} \thispagestyle{plain}
{\begingroup\makeatletter
\let\@makefnmark\relax 
\makeatother\endgroup}
\begin{quote}
\small {\sc Abstract.} In this paper, we show that for exact
area-preserving twist maps on annulus, the invariant circles with a
given rotation number can be destroyed by arbitrarily small
Gevrey-$\alpha$ perturbations of the integrable generating function
in the $C^r$ topology with $r<4-\frac{2}{\alpha}$, where $\alpha>1$.
\end{quote}
\begin{quote}
\end{quote}
\begin{quote}
 \end{quote} \vspace{20pt}

\section{\sc Introduction and main result}
For exact area-preserving twist maps on annulus, it was proved by
Herman in \cite{H2} that invariant circles with given rotation
numbers can be destroyed by $C^{3-\epsilon}$ arbitrarily small
$C^\infty$ perturbations. Following the ideas and techniques
developed by J.N.Mather in the series of papers \cite{M1,M2,M3,M4}, a variational proof of Herman's result was provided in
\cite{W}.  In contrast with it, it has been shown that KAM invariant circles with certain rotation number persist under arbitrarily small
perturbations in the $C^3$ topology (\cite{H33}). For Hamiltonian systems with multi-degrees of freedom, the corresponding results were obtained by \cite{C,CW} and \cite{P}. A partial result on destruction of all invariant tori can be found in \cite{W2}.

On the other hand, for certain rotation numbers, it was obtained by
Mather (resp. Forni) in \cite{M4} (resp. \cite{F}) that the
invariant circles with that rotation numbers can be destroyed by
small perturbations in finer topology respectively. More precisely,
Mather considered Liouville rotation numbers and the topology of the
perturbation induced by $C^\infty$ metric. Forni was concerned about
more special rotation numbers which can be approximated  by rational
ones exponentially fast and the topology of the perturbation induced by
the supremum norm of real-analytic function. Roughly speaking, there
is a balance among the arithmetic property of the rotation number,
the regularity of the perturbation and its topology.

Comparing the results on both sides, it is natural to ask what
happens for perturbations of regularity between $C^\infty$ and
$C^\omega$ (real-analytic). Gevrey-$\alpha$ ($\alpha\geq 1$) functions (see
Definition \ref{gev}) characterize that kind of regularity
quantitatively.  Gevrey Hamiltonians were considered in lots of works (see \cite{MS,MS2} and \cite{pop} for instance). Gevrey-1 functions correspond to ``the best"
$C^\infty$ functions, i.e. $C^\omega$ functions. Gevrey-$\infty$
functions are equivalent to ``the worst" $C^\infty$ functions. For $\alpha>1$, there are compactly supported functions in the class that are not identically zero. This gives much more flexibility to construct examples and show non-existence of invariant circles. In this paper,  we consider the following problem: \begin{itemize}
\item for every given
rotation number $\omega$ and $\alpha$ ($\alpha>1$), what is
the maximum value of $r$ such that the
invariant circle with $\omega$ can be destroyed by an arbitrarily
small Gevrey-$\alpha$ perturbations of the integrable generating function in the $C^r$ topology?
\end{itemize}

To state our result,  we first introduce some terminology. An
irrational number $\omega\in\R$ is called $\mu$-approximated if
there exists a positive number $C>0$ as well as infinitely many
integers $p_n\in\Z$ and $q_n\in\N$ such that
\begin{equation}\label{mu app}
|q_n\omega-p_n|<Cq_n^{-1-\mu}.
\end{equation}
It follows from Dirichlet approximation that any irrational number
is 0-approximated. In particular, $\omega$ is called Liouville if it
is $\mu$-approximated for all $\mu>0$. Given a completely integrable system with
the generating function
\[h_0(x,x')=\frac{1}{2}(x-x')^2, \quad x,x'\in \R,\] we solve the problem above partially. More precisely, we have the
following  theorem.
\begin{Theorem}\label{liou}
For exact area-preserving twist maps on annulus, the invariant
circles with a given $\mu$-approximated rotation number can be
destroyed by arbitrarily small  Gevrey-$\alpha$ ($\alpha> 1$) perturbations of $h_0$ in the $C^r$ topology with
$r<2+\left(2-\frac{2}{\alpha}\right)(1+\mu)$. In particular, the invariant circles with a given Liouville rotation
number can be destroyed by arbitrarily small Gevrey-$\alpha$ ($\alpha> 1$)
perturbations of $h_0$ in the
$C^\infty$ topology.
\end{Theorem}

 Obviously, Theorem \ref{liou} implies Herman's result (\cite{H2}) and Mather's result (\cite{M4}). Unfortunately, we still don't know
whether our result is optimal in the class of  Gevrey-$\alpha$ ($\alpha\geq 1$)
perturbations. Some further developments of
KAM theory are needed to verify the optimality.

For the proof of Theorem \ref{liou}, our approach is parallel to an
investigation  of  variational destruction of invariant circles
under $C^{4-\delta}$ arbitrarily small $C^\infty$ perturbations of
generating functions in \cite{W} (see also \cite{F,M4}). Hence, some parts of the
respective exposition are quite similar. But we decided to repeat
them anyway such that the reader needs not refer to \cite{F,M4,W} for the essentials.
\section{\sc Preliminaries}

\subsection{Minimal configuration}
Let $F$ be a diffeomorphism of $\R^2$ denoted by $F(x,y)=(X(x,y),Y(x,y))$. Let $F$ satisfy:
\begin{itemize}
\item {\it Periodicity:} $F\circ T=T\circ F$ for the translation $T(x,y)=(x+1,y)$;
\item {\it Twist condition:} the map $\psi:(x,y)\mapsto(x,X(x,y))$ is a diffeomorphism of $\R^2$;
\item {\it Exact symplectic:} there exists a real valued function $h$ on $\R^2$ with $h(x+1,y)=h(x,y)$ such that
    \[YdX-ydx=dh.\]
\end{itemize}
Then $F$ induces a map on the cylinder denoted by $f$: $\T\times\R\mapsto \T\times\R$ ($\T=\R/\Z$). $f$ is called an exact
area-preserving monotone twist map. The function  $h$: $\R^2\rightarrow\R^2$ is called
a generating function of $F$, namely $F$
is generated by the following equations
\begin{equation*}
\begin{cases}
y=-\partial_1 h(x,x'),\\
y'=\partial_2 h(x,x'),
\end{cases}
\end{equation*}
where $F(x,y)=(x',y')$.

The function $F$ gives rise to a dynamical
system whose orbits are given by the images of points of $\R^2$
under the successive iterates of $F$. The orbit of the point
$(x_0,y_0)$ is the bi-infinite sequence
\[\{...,(x_{-k},y_{-k}),...,(x_{-1},y_{-1}),(x_0,y_0),(x_1,y_1),...,(x_k,y_k),...\},\]
where $(x_k,y_k)=F(x_{k-1},y_{k-1})$. The sequence
\[(...,x_{-k},...,x_{-1},x_0,x_1,...,x_k,...)\] denoted by $(x_i)_{i\in\Z}$ is called a
stationary configuration if it stratifies the identity
\[\partial_1 h(x_i,x_{i+1})+\partial_2 h(x_{i-1},x_i)=0,\ \text{for\ every\ }i\in\Z.\]
Given a sequence of points $(z_i,...,z_j)$, we can associate its
action
\[h(z_i,...,z_j)=\sum_{i\leq s<j}h(z_s,z_{s+1}).\] A configuration $(x_i)_{i\in\Z}$
is called minimal if for any $i<j\in \Z$, the segment
$(x_i,...,x_j)$ minimizes $h(z_i,...,z_j)$ among all segments
$(z_i,...,z_j)$ of the configuration  satisfying $z_i=x_i$ and
$z_j=x_j$. It is easy to see that every minimal configuration is a
stationary configuration. There is a visual way to describe  configurations. A configuration $(x_i)_{i\in\Z}$ is a function from $\Z$ to $\R$. One can interpolate this function linearly and obtain a piecewise affine function $\R\rightarrow\R$ denoted by $t\mapsto x_t$. The graph of this function is sometimes called the Aubry diagram of the configuration.
By \cite{B} (see also \cite{G}), minimal configurations satisfy a
group of remarkable properties as follows:
\begin{itemize}
\item Two distinct minimal configurations seen as the Aubry diagrams cross at most once, which
is so called Aubry's crossing lemma.
\item For every minimal configuration $\bold{x}=(x_i)_{i\in\Z}$, the limit
\[\rho(\bold{x})=\lim_{n\rightarrow\infty}\frac{x_{i+n}-x_i}{n}\]
exists and doesn't depend on $i\in\Z$. $\rho(\bold{x})$ is called
the rotation number of $\bold{x}$.
\item For every $\omega\in \R$, there exists a minimal configuration
with rotation number $\omega$. Following the notations of \cite{B}, the
set of all minimal configurations with rotation number $\omega$ is
denoted by $M_\omega^h$, which can be endowed with the topology
induced from the product topology on $\R^\Z$. If
$\bold{x}=(x_i)_{i\in\Z}$ is a minimal configuration, considering
the projection $pr:\ M_\omega^h\rightarrow\R$ defined by
$pr(\bold{x})=x_0$, we set $\mathcal {A}_\omega^h=pr(M_\omega^h)$.
\item If $\omega\in\Q$, say $\omega=p/q$ (in lowest terms), then it is convenient to define the rotation symbol to detect the structure of
$M_{p/q}^h$. If $\bold{x}$ is a minimal configuration with rotation
number $p/q $, then the rotation symbol $\sigma(\bold{x})$ of
$\bold{x}$ is defined as follows
\begin{equation*}
\sigma(\bold{x})=\left\{\begin{array}{ll}
\hspace{-0.4em}p/q+,&\text{if}\ x_{i+q}>x_i+p\ \text{for\ all\ }i,\\
\hspace{-0.4em}p/q,&\text{if}\ x_{i+q}=x_i+p\ \text{for\ all\ }i,\\
\hspace{-0.4em}p/q-,&\text{if}\ x_{i+q}<x_i+p\ \text{for\ all\ }i.\\
\end{array}\right.
\end{equation*}
 Moreover, we set
\begin{align*}
&M_{{p/q}^+}^h=\{\bold{x} \text{\  is a minimal configuration with
rotation symbol}\  p/q \text{\ or\ } p/q+\},\\
&M_{{p/q}^-}^h=\{\bold{x} \text{\ is a minimal configuration with
rotation symbol}\  p/q \text{\ or\ } p/q-\},
\end{align*}
then both $M_{{p/q}^+}^h$ and $M_{{p/q}^+}^h$ are totally ordered.
Namely, every two configurations in each of them (seen as Aubry diagrams) do not cross. We
denote $pr(M_{{p/q}^+}^h)$ and $pr(M_{{p/q}^-}^h)$ by $\mathcal
{A}_{{p/q}^+}^h$ and $\mathcal {A}_{{p/q}^-}^h$ respectively.
\item If $\omega\in\R\backslash\Q$ and $\bold{x}$ is a minimal
configuration with rotation number $\omega$, then
$\sigma(\bold{x})=\omega$ and $M_\omega^h$ is totally ordered.
\item $\mathcal {A}_\omega^h$ is a closed subset of $\R$ for every rotation symbol
$\omega$.
\end{itemize}
\subsection{Peierls's barrier}
In \cite{M3}, Mather introduced the notion of Peierls's barrier and gave
a criterion of existence of invariant circle. Namely, the exact
area-preserving monotone twist map generated by $h$ admits an
invariant circle with rotation number $\omega$ if and only if the
Peierls's barrier $P_\omega^h(\xi)$ vanishes identically for all
$\xi\in\R$. The Peierls's barrier is defined as follows:
\begin{itemize}
\item If $\xi\in \mathcal {A}_\omega^h$, we set $P_\omega^h(\xi)$=0.
\item If $\xi \not\in \mathcal {A}_\omega^h$, since $\mathcal {A}_\omega^h$ is a closed set in $\R$, then $\xi$ is contained in some
complementary interval $(\xi^-,\xi^+)$ of $\mathcal {A}_\omega^h$ in
$\R$. By the definition of $\mathcal {A}_\omega^h$, there exist
minimal configurations with rotation symbol $\omega$,
$\bold{x^-}=(x_i^-)_{i\in\Z}$ and $\bold{x^+}=(x_i^+)_{i\in\Z}$
satisfying $x_0^-=\xi^-$ and $x_0^+=\xi^+$. For every configuration
$\bold{x}=(x_i)_{i\in\Z}$ satisfying $x_i^-\leq x_i\leq x_i^+$, we
set
\[G_\omega(\bold{x})=\sum_I(h(x_i,x_{i+1})-h(x_i^-,x_{i+1}^-)),\]
where $I=\Z$, if $\omega$ is not a rational number, and $I=\{0,...,
q-1 \}$, if $\omega=p/q$. $P_\omega^h(\xi)$ is defined as the
minimum of $G_\omega(\bold{x})$ over the configurations $\bold{x}\in
\Pi=\prod_{i\in I}[x_i^-,x_i^+]$ satisfying $x_0=\xi$. Namely
\[P_\omega^h(\xi)=\min_{\bold{x}}\{G_\omega(\bold{x})|\bold{x}\in \Pi\ \text{and}\ \ x_0=\xi\}.\]
\end{itemize}
By \cite{M3}, $P_\omega^h(\xi)$ is a non-negative periodic function of
the variable $\xi\in\R$ with the modulus of continuity with respect
to $\omega$ and its modulus of continuity with respect to $\omega$ can be bounded from above. Due to the periodicity of $P_\omega^h(\xi)$ with
respect to $\xi$, we only need to consider it in the interval $[0,1]$.

\subsection{Gevrey function}
We fix $\alpha>1$ and let $K$ be a closed interval in $\R$. We
recall the definitions
\begin{equation}
\|\phi\|_{\alpha,L}:=\sum_{k\in\N}\frac{L^{|k|\alpha}}{k!^\alpha}\|\partial^k\phi\|_{C^0(K)}
\end{equation}
\[G^{\alpha,L}(K)=\{\phi\in C^\infty(K)| \|\phi\|_{\alpha,L}<\infty\},\quad G^\alpha(K)=\bigcup_{L>0}G^{\alpha,L}(K).\]
Following \cite{MS}, Gevrey-$\alpha$ function is defined as follow.
\begin{Definition}\label{gev}
A function $\phi$ is called Gevrey-$\alpha$ function on $K$ if  $\phi\in G^\alpha(K)$.
\end{Definition}
 From Leibniz rule, it follows that for $L>0$, $\phi,\psi\in
G^{\alpha,L}(K)$,
\begin{equation}\label{leib}
\|\phi\psi\|_{\alpha,L}\leq \|\phi\|_{\alpha,L}\|\psi\|_{\alpha,L}.
\end{equation}

 For the simplicity of notations, we don't distinguish the constant $C$ in
following different estimate formulas.
\section{\sc Construction of the generating functions}
In order to destroy  the
invariant circle with  a given rotation number of the completely integrable system
\[h_0(x,x')=\frac{1}{2}(x-x')^2, \quad x,x'\in \R,\] we construct the perturbation consisting of two
parts. The first one is
\begin{equation}\label{31}
u_n(x)=\frac{1}{n^a}(1-\cos(2\pi x) ),\quad x\in \R,\end{equation}
where $n\in \N$ and $a$ is a positive constant independent of $n$.

We construct  the second part of the perturbation in the following.
First of all, for each $\lambda>0$, we construct a function
$f_\lambda\in C^\infty(\R)$ as follow:
\begin{equation}
f_\lambda(x)=\left\{\begin{array}{ll}\hspace{-0.4em}0,&x\leq 0,\\
\hspace{-0.4em}\exp\left(-\lambda\sqrt{2}x^{-\frac{1}{\alpha-1}}\right),&\text{otherwise}.\\
\end{array}\right.
\end{equation}
\begin{Lemma}
There exists $\lambda>0$ such that $f_\lambda(x)$ is a
Gevrey-$\alpha$ function on $\R$.
\end{Lemma}
\Proof For the simplicity of notations, let $p=\frac{1}{\alpha-1}$.
Since $\alpha\in (1,\infty)$, $p\in (0,\infty)$ and
\begin{equation}
f_\lambda(x)=\left\{\begin{array}{ll}\hspace{-0.4em}0,&x\leq 0,\\
\hspace{-0.4em}\exp\left(-\lambda\sqrt{2}x^{-p}\right),&\text{otherwise}.\\
\end{array}\right.
\end{equation}

Let $k\in\N$ and $x>0$. We observe that $f_\lambda|\R^+$ can be
extended to a holomorphic function on $\C\backslash (-\infty,0]$. Let
$\sigma:=\frac{\pi}{4}\min\{1,\frac{1}{p}\}$ and
$\Sigma_\sigma=\{z\in\C||\text{arg}z|\leq \sigma\}$. The closed disk
$D_z$ of center $x$ and radius $(x\sin\sigma)$ is the largest disk
centered at $x$ contained in $\Sigma_\sigma$, and the Cauchy
inequalities yield
\[\left|f_\lambda^{(k)}(x)\right|\leq\frac{k!}{(x\sin\sigma)^k}\max_{D^z}|f_\lambda|.\]
Let $z=re^{i\theta}\in D_z$. Since $|\theta|\leq \sigma=\frac{\pi}{4}\min\{1,\frac{1}{p}\}$, then
$\mathfrak{R}e(z^{-p})=r^{-p}\cos(p\theta)\geq
\frac{1}{\sqrt{2}|z|^p}$ and $|z|\leq 2x$. Hence, we have
\[\max_{D_z}|f_\lambda|\leq \exp\left(-\frac{\lambda}{(2x)^p}\right).\]
It is easy to see that the maximum of $y\mapsto y^ke^{-\lambda
y^p}$ is $\left(\frac{k}{\lambda p e}\right)^{k/p}$, therefore
\begin{equation}\label{kd}
\left|f_\lambda^{(k)}(x)\right|\leq\left(\frac{2}{\sin\sigma}\right)^k\left(\frac{k}{\lambda p e}\right)^{k/p}k!.
\end{equation}
By Stirling formula, we have that for any given $L>0$, if
$\lambda>(2L^\alpha/\sin\sigma)^p/p$, then
\[\sum_{k\in\N}\frac{L^{|k|\alpha}}{k!^\alpha}\left\|f_\lambda^{(k)}(x)\right\|_{C^0(\R)}<\infty.\]
Therefore, there exists $\lambda>0$ such that $f_\lambda(x)$ is a
Gevrey-$\alpha$ function on $\R$.\End

$v_n(x)$ is constructed as follows. For $x\in [0,1]$, we let
\begin{equation*}
v_n(x)=\left\{\begin{array}{ll}\hspace{-0.4em}\frac{1}{n^a}f_\lambda\left(\frac{1}{8n^{a/2}}-\frac{1}{2}+x\right)f_\lambda
\left(\frac{1}{8n^{a/2}}+\frac{1}{2}-x\right),&x\in \left[\frac{1}{2}-\frac{1}{8n^{a/2}},\frac{1}{2}+\frac{1}{8n^{a/2}}\right],\\
\hspace{-0.4em}0,&[0,1]\backslash\left[\frac{1}{2}-\frac{1}{8n^{a/2}},\frac{1}{2}+\frac{1}{8n^{a/2}}\right].\\
\end{array}\right.
\end{equation*}
More precisely, for $x\in
\left[\frac{1}{2}-\frac{1}{8n^{a/2}},\frac{1}{2}+\frac{1}{8n^{a/2}}\right]$,
\begin{equation*}
v_n(x)=\frac{1}{n^a}\exp\left(-\lambda\sqrt{2}\left(\left(\frac{1}{8n^{a/2}}-\frac{1}{2}+x\right)^{-\frac{1}{\alpha-1}}+\left(\frac{1}{8n^{a/2}}+
\frac{1}{2}-x\right)^{-\frac{1}{\alpha-1}}\right)\right),
\end{equation*}
where $\lambda$ is a positive constant independent of $n$. Moreover,
we extend $v_n(x)$ on $[0,1]$ to be a periodic function on $\R$ by
$v_n(x+1)=v_n(x)$. By (\ref{leib}), $v_n(x)$ is  a Gevrey-$\alpha$
function on $\R$. Based on the definition of $v_n$, it follows from a simple
calculation that for $\alpha\in (1,\infty)$
\begin{equation}\label{maxv}
\max_{[0,1]}v_n(x)=v_n\left(\frac{1}{2}\right)\sim\frac{1}{n^a}\exp\left(-Cn^{\frac{a}{2(\alpha-1)}}\right),
\end{equation}
where  $f\sim g$ means that $\frac{1}{C}g<f<C g$ holds for some
constant $C>0$. From (\ref{kd}), it follows that for $r>0$, we have
\begin{equation}\label{ff}
\left\|f_\lambda^{(r)}\left(\frac{1}{8n^{a/2}}-\frac{1}{2}+x\right)\right\|_{C^0}\leq C_1 \ \text{and}\ \left\|f_\lambda^{(r)}\left(\frac{1}{8n^{a/2}}+\frac{1}{2}-x\right)\right\|_{C^0}\leq C_1,
\end{equation}
where $C_1$ is a positive constant independent of $n$. By Leibniz formula, we have
\begin{equation*}
\|v_n^{(r)}(x)\|_{C^0}\leq \frac{1}{n^a}\sum_{i=0}^rC_r^i\left\|f_\lambda^{(i)}\left(\frac{1}{8n^{a/2}}-\frac{1}{2}+x\right)\right\|_{C^0}\cdot
\left\|f_\lambda^{(r-i)}
\left(\frac{1}{8n^{a/2}}+\frac{1}{2}-x\right)\right\|_{C^0} ,
\end{equation*}
which together with (\ref{ff}) implies that for any fixed $r>0$ and $n$ large enough, we have
\begin{equation}\label{vr}
\|v_n(x)\|_{C^r([0,1])}\leq C_2\frac{1}{n^a},
\end{equation}
where $C_2$ is a positive constant independent of $n$.

So far, we complete the construction of the generating function of
the nearly integrable system,
\begin{equation}\label{h}
h_n(x,x')=h_0(x,x')+u_n(x')+v_n(x'),
\end{equation}
where $n\in\N$.

\section{\sc Proof of Theorem \ref{liou}}
 If $\omega\in \Q$, then the invariant circles with rotation
number $\omega$ could be easily destroyed by an analytic
perturbation arbitrarily close to zero. Therefore it suffices to
consider the irrational $\omega$. Firstly, we prove the
non-existence of invariant circles with a small enough rotation
number. More precisely, we have the following Lemma:
\begin{Lemma}\label{MR} For $\omega\in\R\backslash\Q$ and $n$ large enough,
the exact area-preserving monotone twist map generated by $h_n$
admits no invariant circle with the rotation number satisfying
\[|\omega|<n^{-\frac{a\alpha}{2(\alpha-1)}-\delta}, \] where $\delta$ is a small positive constant independent of $n$.
\end{Lemma}
First of all, we will estimate the lower bound of $P_{0^+}^{h_n}$ at
a given point. To achieve that, we need to estimate the distances of
pairwise adjacent elements of the minimal configuration. More
precisely, we have
\begin{Lemma}\label{lowstep} Let $(x_i)_{i\in \Z}$ be a minimal
configuration of $\bar{h}_n$ with rotation symbol $\omega>0$, then
\[x_{i+1}-x_i\geq \frac{1}{2}n^{-\frac{a}{2}},\quad \text{for}\quad x_i\in \left[\frac{1}{4},\frac{3}{4}\right].\]\end{Lemma}

\Proof Without loss of generality, we assume $x_i\in [0,1]$ for all
$i\in\Z$. By Aubry's crossing lemma, we have
\[...<x_{i-1}<x_i<x_{i+1}<....\]We consider the configuration
$(\xi_i)_{i\in \Z}$ defined by
\begin{equation*}
\xi_j= \left\{\begin{array}{ll}\hspace{-0.4em}x_j,& j<i,\\
\hspace{-0.4em}x_{j+1},& j\geq i.\\
\end{array}\right.
\end{equation*}
Since $(x_i)_{i\in \Z}$ is minimal, we have
\[\sum_{i\in \Z}\bar{h}_n(\xi_i,\xi_{i+1})-\sum_{i\in \Z}\bar{h}_n(x_i,x_{i+1})\geq 0.\]
By the definitions of $\bar{h}_n$ and $(\xi_i)_{i\in\Z}$, we have
\begin{align*}
0&\leq\sum_{i\in \Z}\bar{h}_n(\xi_i,\xi_{i+1})-\sum_{i\in
\Z}\bar{h}_n(x_i,x_{i+1})\\
&=\bar{h}_n(x_{i-1},x_{i+1})-\bar{h}_n(x_{i-1},x_{i})-\bar{h}_n(x_{i},x_{i+1})\\
&=(x_{i+1}-x_i)(x_i-x_{i-1})-u_n(x_i).
\end{align*}
Moreover,\[u_n(x_i)\leq(x_{i+1}-x_i)(x_i-x_{i-1})\leq\frac{1}{4}(x_{i+1}-x_{i-1})^2.\]
Therefore,
\begin{equation}\label{uwith}
x_{i+1}-x_{i-1}\geq 2\sqrt{u_n(x_i)}. \end{equation} For $x_i\in
[\frac{1}{4},\frac{3}{4}]$, $u_n(x_i)\geq n^{-a}$, hence,
\begin{equation}\label{ls} x_{i+1}-x_{i-1}\geq 2n^{-\frac{a}{2}}.
\end{equation}

Since $(x_i)_{i\in\Z}$ is a stationary configuration, we have
\begin{align*}
x_{i+1}-x_i&=-\partial_1\bar{h}_n(x_i,x_{i+1}),\\
&=\partial_2\bar{h}_n(x_{i-1},x_i),\\
&=x_i-x_{i-1}+u_n'(x_i).
\end{align*}
Since $u_n'(x)=\frac{2\pi}{n^a}\sin(2\pi x)$, it follows from
$(\ref{ls})$ that
\[x_{i+1}-x_i\geq \frac{1}{2}n^{-\frac{a}{2}},\quad x_i\in \left[\frac{1}{4},\frac{3}{4}\right]
.\] The proof of Lemma \ref{lowstep} is completed.\End

In order to estimate the lower bound of $P_{0^+}^{h_n}$ at a given point, we make  a modification of  $v_n$ by changing the axis of symmetry of its support into $\eta$. We denote $v_{n,\eta}(x):=v_n(x-(\eta-1/2))$. Then we have
\[\text{supp}\ v_{n,\eta}=\left[\eta-\frac{1}{8n^{a/2}},\eta+\frac{1}{8n^{a/2}}\right]\]
 By a similar calculation as (\ref{maxv}), we have
\begin{equation}
v_{n,\eta}\left(\eta\right)=\max
v_{n,\eta}(x)\sim\frac{1}{n^a}\exp\left(-Cn^{\frac{a}{2(\alpha-1)}}\right).
\end{equation}

 Let $(x_i)_{i\in \Z}$ be the minimal
configuration of
$\bar{h}_n(x_i,x_{i+1})=h_0(x_i,x_{i+1})+u_n(x_{i+1})$ with rotation
symbol $0^+$, then
from Lemma \ref{lowstep}, we have
\[x_{i+1}-x_i\geq \frac{1}{2}n^{-\frac{a}{2}},\quad x_i\in \left[\frac{1}{4},\frac{3}{4}\right]
.\]
Hence, there exists $\eta\in \left[\frac{3}{8},\frac{5}{8}\right]$ such that
\[(x_i)_{i\in \Z}\cap\text{supp}\ v_{n,\eta}=\emptyset.\]
Moreover, for all $i\in\Z$, \[v_{n,\eta}(x_i)=0.\]

Based on \cite{M4} (p.207-208), the Peierls's barrier $P_{0^+}^{h_n}(\eta)$
could be defined as follows
\begin{equation*}
P_{0^+}^{h_n}(\eta)=\min_{\xi_0=\eta}\sum_{i\in
\Z}h_n(\xi_i,\xi_{i+1})-\min\sum_{i\in \Z}h_n(z_i,z_{i+1}),
\end{equation*}
where $(\xi_i)_{i\in\Z}$ and  $(z_i)_{i\in\Z}$ are monotone
increasing configurations limiting on $0,\ 1$.
 Let $(\xi_i)_{i\in \Z}$ and  $(z_i)_{i\in \Z}$
be  minimal configurations of $h_n$ defined by (\ref{h}) with
rotation symbol $0^+$ satisfying $\xi_0=\eta$ and . Then we have
\begin{align*}
\sum_{i\in
\Z}(h_n(&\xi_i,\xi_{i+1})-h_n(z_i,z_{i+1}))\\
&\geq v_{n,\eta}\left(\eta\right)+\sum_{i\in
\Z}\bar{h}_n(\xi_i,\xi_{i+1})-\sum_{i\in
\Z}h_n(z_i,z_{i+1}),\\
&\geq v_{n,\eta}\left(\eta\right)+\sum_{i\in
\Z}\bar{h}_n(x_i,x_{i+1})-\sum_{i\in
\Z}h_n(z_i,z_{i+1}),\\
&\geq v_{n,\eta}\left(\eta\right)+\sum_{i\in
\Z}\bar{h}_n(x_i,x_{i+1})
-\sum_{i\in
\Z}h_n(x_i,x_{i+1}),\\
&=v_{n,\eta}\left(\eta\right)-\sum_{i\in \Z}v_{n,\eta}(x_{i+1}),\\
&=v_{n,\eta}\left(\eta\right),
\end{align*}
where the first inequality holds since $v_{n,\eta}\geq 0$, the second one since $(x_i)_{i\in \Z}$ is a minimal
configuration of
$\bar{h}_n$, the third one  since $(z_i)_{i\in \Z}$
is a minimal configuration of $h_n$ and the last one since $v_{n,\eta}(x_i)=0$ for all $i\in\Z$.
Moreover, we have
\[P_{0^+}^{h_n}\left(\eta\right)\geq
v_{n,\eta}\left(\eta\right).\] It follows that
\begin{equation}\label{lowb}
P_{0^+}^{h_n}\left(\eta\right)\geq
\frac{C_1}{n^a}\exp\left(-C_2n^{\frac{a}{2(\alpha-1)}}\right).
\end{equation}

 Second, following a similar argument as \cite{W}, one can obtain the improvement of modulus of continuity of Peierls's barrier based on the hyperbolicity of $h_n$. More precisely, we have the following lemma.
 \begin{Lemma}\label{pw}
 For every irrational rotation symbol $\omega$ satisfying $0<\omega<n^{-\frac{a\alpha}{2(\alpha-1)}-\delta}$, we have
\begin{equation}\label{app}
\left|P_{\omega}^{h_n}\left(\eta\right)-P_{0^+}^{h_n}\left(\eta\right)\right|\leq
C_1\exp\left(-C_2n^{\frac{a}{2(\alpha-1)}+\frac{\delta}{2}}\right).\end{equation}where $\eta\in [3/8,5/8]$ and
$\delta$ is a small positive constant independent of $n$.
\end{Lemma}
\Proof
 If $\eta\in
\mathcal {A}_\omega^{h_n}$, then $P_\omega^{h_n}(\eta)=0$. Hence, it
suffices to consider the case with $\eta\not\in \mathcal
{A}_\omega^{h_n}$ to destroy invariant circles. Since the proof of Lemma \ref{pw} is similar to Lemma 5.1 in \cite{W}, we only give a sketch of the proof to show some main differences between them. For the simplicity of notations, we denote $\kappa:=\frac{a}{2(\alpha-1)}$ and $\epsilon_n:=\exp(-n^{\kappa+\frac{\delta}{2}})$. The proof is proceeded by  three steps as follows.

In the first step, we will show that each of the intervals $[0,\epsilon_n]$ and $[1-\epsilon_n,1]$ contains a large number of elements of the minimal configuration $(x_i)_{i\in\Z}$ of $h_n$ with irrational rotation symbol $0<\omega<n^{-\kappa-\frac{a}{2}-\delta}$ for $n$ large enough. Let \[\Sigma_n=\left\{i\in \Z\ |\ \,x_i\in \left[\epsilon_n,1-\epsilon_n\right]\right\},\] then it follows from a similar argument as Lemma 5.2 in \cite{W} that
\begin{equation}\label{sig}
\sharp\Sigma_n\leq
Cn^{\kappa+\frac{a}{2}+\frac{\delta}{2}},
\end{equation}
where $\sharp\Sigma_n$ denotes
the number of elements in $\Sigma_n$.  Let $I$ be a interval of length
$1$. We denote $\Delta_\omega:=\{i\in \Z\ |\ x_i\in
I\}$. Since $(x_i)_{i\in \Z}$ is a minimal configuration with rotation number
$\omega\in \R\backslash\Q$, then it follows from  Lemma 5.3 in \cite{W} that
\begin{equation}
\frac{1}{\omega}-1\leq\sharp\Delta_\omega\leq\frac{1}{\omega}+1.
\end{equation} which together with (\ref{sig}) implies
\[\sharp\Delta_\omega\geq Cn^{\kappa+\frac{a}{2}+\delta}\gg Cn^{\kappa+\frac{a}{2}+\frac{\delta}{2}}\geq \sharp\Sigma_n.\]
Hence, one can obtain that each of the intervals $[0,\epsilon_n]$ and $[1-\epsilon_n,1]$ contains a large number of elements of the minimal configuration $(x_i)_{i\in\Z}$ of $h_n$ with irrational rotation symbol $0<\omega<n^{-\kappa-\frac{a}{2}-\delta}$ for $n$ large enough (see Lemma 5.4 in \cite{W}).

In the second step,  we approximate
$P_\omega^{h_n}(\eta)$ for $\eta\in
[3/8,5/8]$ by
the difference of the actions of the segments with a given length, where we consider the number of the elements in a segment of the configuration as the length of the segment. Let
$(\xi^-,\xi^+)$ be the complementary interval of $\mathcal
{A}_\omega^{h_n}$ in $\R$ and contains $\eta$. Let
$(\xi_i^{\pm})_{i\in\Z}$ be the minimal
configurations with rotation symbol $\omega$ satisfying
$\xi_0^{\pm}=\xi^{\pm}$ and  let $(\xi_i)_{i\in \Z}$ be a minimal
configuration of $h_n$ with rotation symbol $\omega$ satisfying
$\xi_0=\eta$ and $\xi_i^-\leq \xi_i\leq \xi_i^+$.
We denote
$d(x):=\min\{|x|,|x-1|\}$. By Step 1, there exist $i^-,\ i^+$ such
that
\begin{equation}\label{5}
d(\xi_i^-)<\epsilon_n\quad \text{and}\quad
\xi_{i+1}^--\xi_{i-1}^-\leq\epsilon_n\quad \text{for}\quad i=i^-,\
i^+. \end{equation}
 Thanks to Aubry's crossing lemma, we have $\xi_i^-\leq
\xi_i\leq \xi_i^+\leq \xi_{i+1}^-$. Hence,
\[\xi_i-\xi_i^-\leq\epsilon_n\quad \text{for}\quad i=i^-,\ i^+.\]
We define the following configuration:
\begin{align*}
y_i=\left\{\begin{array}{ll}
\hspace{-0.4em}\xi_i,& i^-<i<i^+,\\
\hspace{-0.4em}\xi_i^-,& i\leq i^-,\ i\geq i^+.
\end{array}\right.
\end{align*}
Since $\eta\in
[3/8,5/8]\subset [\epsilon_n,1-\epsilon_n]$ for $n$ large enough, then  $\xi_0=\eta$ is contained in  $(y_i)_{i\in\Z}$ up to the rearrangement of the
index $i$. By a direct calculation (see (11)-(15) in \cite{W}), we have
\begin{equation}\label{step1}
P_\omega^{h_n}(\eta)\leq\sum_{i\in
\Z}(h_n(y_i,y_{i+1})-h_n(\xi_i^-,\xi_{i+1}^-))\leq
P_\omega^{h_n}(\eta)+C\epsilon_n^2,
\end{equation}

In the third step, we will compare $P_{0^+}^{h_n}(\eta)$ with $\sum_{i\in
\Z}(h_n(y_i,y_{i+1})-h_n(\xi_i^-,\xi_{i+1}^-))$.
By \cite{M4}, we have
\begin{equation*}
P_{0^+}^{h_n}(\eta)=\min_{\xi_0=\eta}\sum_{i\in
\Z}h_n(\xi_i,\xi_{i+1})-\min\sum_{i\in \Z}h_n(z_i,z_{i+1}),
\end{equation*}
where $(\xi_i)_{i\in\Z}$ and  $(z_i)_{i\in\Z}$ are monotone
increasing configurations limiting on $0,\ 1$. We denote
\begin{equation*}
\begin{cases}
K(\eta)=\min_{\xi_0=\eta}\sum_{i\in \Z}h_n(\xi_i,\xi_{i+1}),\\
K=\min\sum_{i\in \Z}h_n(z_i,z_{i+1}).
\end{cases}
\end{equation*}By a direct calculation (see (17)-(28) in \cite{W}), we have
\begin{equation}\label{13}
\left|\sum_{i\in
\Z}h_n(\xi_i^-,\xi_{i+1}^-)-K\right|\leq C\epsilon_n^2\ \ \text{and}\ \ \left|\sum_{i\in
\Z}h_n(y_i,y_{i+1})-K(\eta)\right|\leq C\epsilon_n^2.
\end{equation}
Finally, from $(\ref{step1})$ and $(\ref{13})$, we obtain
\begin{align*}
|P_{\omega}^{h_n}(\eta)-P_{0^+}^{h_n}(\eta)|&\leq|\sum_{i\in
\Z}h_n(y_i,y_{i+1})-\sum_{i\in
\Z}h_n(\xi_i^-,\xi_{i+1}^-)+K-K(\xi)|+C_1\epsilon_n^2,\\
&\leq
|\sum_{i\in
\Z}h_n(y_i,y_{i+1})-K(\xi)|+|\sum_{i\in
\Z}h_n(\xi_i^-,\xi_{i+1}^-)-K|+C_1\epsilon_n^2,\\
&\leq C\epsilon_n^2.
\end{align*}
Recalling $\kappa:=\frac{a}{2(\alpha-1)}$ and $\epsilon_n:=\exp(-n^{\kappa+\frac{\delta}{2}})$, we have
\begin{equation}
\left|P_{\omega}^{h_n}\left(\eta\right)-P_{0^+}^{h_n}\left(\eta\right)\right|\leq
C_1\exp\left(-C_2n^{\frac{a}{2(\alpha-1)}+\frac{\delta}{2}}\right).\end{equation}
which completes the proof of Lemma \ref{pw}.
\End

 Based on
the preparations above, it is easy to prove Lemma \ref{MR}. We
assume that there exists an invariant circle with rotation number
$0<\omega<n^{-\frac{a\alpha}{2(\alpha-1)}-\delta}$ for $h_n$, then
$P_\omega^{h_n}(\xi)\equiv 0$ for every $\xi\in \R$. By Lemma \ref{pw},
we have
\begin{equation}\label{pp}
\left|P_{0^+}^{h_n}\left(\eta\right)\right|\leq
C_1\exp\left(-C_2n^{\frac{a}{2(\alpha-1)}+\frac{\delta}{2}}\right).
\end{equation}

On the other hand, $(\ref{lowb})$ implies that
\[P_{0^+}^{h_n}\left(\eta\right)\geq
\frac{C_1}{n^a}\exp\left(-C_2n^{\frac{a}{2(\alpha-1)}}\right).\]Hence, we have
\[\frac{C_1}{n^a}\exp\left(-C_1n^{\frac{a}{2(\alpha-1)}}\right)\leq C_2\exp\left(-C_2n^{\frac{a}{2(\alpha-1)}+\frac{\delta}{2}}\right).\] It is an obvious
contradiction for $n$ large enough. Therefore, there exists no
invariant circle with rotation number
$0<\omega<n^{-\frac{a\alpha}{2(\alpha-1)}-\delta}$.

For $-n^{-\frac{a\alpha}{2(\alpha-1)}-\delta}<\omega<0$, by
comparing $P_{\omega}^{h_n}(\xi)$ with $P_{0^-}^{h_n}(\xi)$, the
proof is similar. We omit the details. Therefore, the proof of
Theorem \ref{MR} is completed.\End

The case with a given irrational rotation number can be easily
reduced to the one with a small enough rotation number. More
precisely,

\begin{Lemma}\label{Herm} Let $h_P$ be a generating function as follow
\[h_P(x,x')=h_0(x,x')+P(x'),\] where $P$ is a periodic
function of periodic $1$. Let $Q(x)=q^{-2}P(qx),q\in \N$, then the
exact area-preserving monotone twist map generated by
$h_Q(x,x')=h_0(x,x')+Q(x')$ admits an invariant circle with rotation
number $\omega \in \R\backslash \Q$ if and only if the exact
area-preserving monotone twist map generated by $h_P$ admits an
invariant circle with rotation number $q\omega-p, p\in \Z$.
\end{Lemma}

We omit the proof and for more details, see \cite{H2}. For the sake of
simplicity of notations, we denote $Q_{q_n}$ by $Q_n$ and the same
to $u_{q_n}, v_{q_n}$ and $h_{q_n}$. Let
\[Q_n(x)={q_n}^{-2}(u_n(q_nx)+v_n(q_nx)),\] where $(q_n)_{n\in \N}$
is a sequence satisfying (\ref{mu app})
\begin{equation}\label{diri}
|q_n\omega-p_n|<\frac{C}{q^{1+\mu}_n}, \end{equation} where $p_n\in \Z$ and
$q_n\in \N$. Since $\omega\in\R\backslash\Q$, we say
$q_n\rightarrow\infty$ as $n\rightarrow\infty$. Let
$\tilde{h}_n(x,x')=h_0(x,x')+Q_n(x')$, we prove Theorem \ref{liou}
for $(\tilde{h}_n)_{n\in\N}$ as follow:

\Proof Based on Lemma \ref{MR} and
(\ref{diri}), it suffices to take
\[\frac{C}{q^{1+\mu}_n}\leq\frac{1}{{q_n}^{\frac{a\alpha}{2(\alpha-1)}+\delta}},\]
which implies
\begin{equation}\label{33}
a\leq \left(2-\frac{2}{\alpha}\right)(1+\mu)-\epsilon,
\end{equation} where
$\epsilon=2\delta(\frac{\alpha-1}{\alpha})$ and $\delta$ is a small
positive constant independent of $n$. From the constructions of
$u_n$ and $v_n$, it follows from (\ref{31}) and (\ref{vr}) that
\begin{align*}
 ||\tilde{h}_n&(x,x')-h_0(x,x')||_{C^r}\\
 &=||Q_n(x')||_{C^r},\\
&\leq{q_n}^{-2}(||u_n(q_nx')||_{C^r}+||v_n(q_nx')||_{C^r}),\\
&\leq{q_n}^{-2}({q_n}^{-a}(2\pi)^r{q_n}^r+C_1{q_n}^{-a}{q_n}^r),\\
&\leq C_2{q_n}^{r-a-2},
\end{align*}
where $C_1, C_2$ are positive constants only depending on $r$.

To complete the proof, it is enough to make $r-a-2<0$, which
together with (\ref{33}) implies
\[r<a+2\leq 2+\left(2-\frac{2}{\alpha}\right)(1+\mu)-\epsilon.\]This completes the
proof of Theorem \ref{liou} if we take
$a=\left(2-\frac{2}{\alpha}\right)(1+\mu)-\epsilon$ and
$r=2+\left(2-\frac{2}{\alpha}\right)(1+\mu)-2\epsilon$.\End

 \vspace{2ex}
\noindent\textbf{Acknowledgement} The author sincerely
thanks the referees for their careful reading of the manuscript and
invaluable comments which were very helpful in improving this paper. The author also would like to thank
Prof. C.-Q. Cheng for many helpful discussions. This work is under
the support of  the NNSF of China (Grant No. 11171071, 
11171146).

{\sc School of Mathematical Sciences, Fudan University,
Shanghai 200433,
China.}

 {\it E-mail address:} \texttt{linwang.math@gmail.com}

\end{document}